\documentclass[11pt]{article}

\usepackage{graphicx,subfigure}
\usepackage{epsfig}
\usepackage{amssymb}
\usepackage{mathrsfs}
\usepackage{color}
\usepackage[dvipdfm]{hyperref}
\usepackage{amsmath}
\usepackage{ifpdf}

\newcommand{\E}{{\cal E}}

\newtheorem{theorem}{Theorem}[section]
\newtheorem{definition}{Definition}[section]
\newtheorem{lemma}[theorem]{Lemma}
\newtheorem{corollary}[theorem]{Corollary}

\def\whitebox{{\hbox{\hskip 1pt
 \vrule height 6pt depth 1.5pt
 \lower 1.5pt\vbox to 7.5pt{\hrule width
    3.2pt\vfill\hrule width 3.2pt}%
 \vrule height 6pt depth 1.5pt
 \hskip 1pt } }}
\def\qed{\ifhmode\allowbreak\else\nobreak\fi\hfill\quad\nobreak
     \whitebox\medbreak}

\newcommand{\ignore}[1]{}

\begin {document}

\baselineskip 16pt
\title{Further results
regarding  the degree resistance distance of cacti}

 \author{\small   Jia\textrm{-}Bao \ Liu$^{a,b}$, \ \ Wen\textrm{-}Rui Wang$^{a}$, \ \ Yong\textrm{-}Ming Zhang$^{a}$, \ \ Xiang\textrm{-}Feng \ Pan$^{a,}$\thanks{Corresponding author. Tel:+86-551-63861313.
  \E-mail:liujiabaoad@163.com(J.Liu), Ricciawang@163.com(W.Wang), ymzhang625@163.com(Y.Zhang), xfpan@ahu.edu.cn(X.Pan).}\\
 \small  $^{a}$ School of Mathematical Sciences, Anhui  University, Hefei 230601, P. R. China\\
 \small  $^{b}$ Department of Public Courses, Anhui Xinhua
 University, Hefei 230088, P. R. China\\}

\date{}
\maketitle
\begin{abstract}
  A graph $G$ is called a cactus if each block of $G$ is either an edge or
a cycle. Denote by $Cact(n;t)$ the set of connected cacti
possessing $n$ vertices and $t$ cycles. In this paper, we show
that there are some errors in [J. Du, G. Su, J. Tu, I. Gutman, The
degree resistance distance of cacti, Discrete Appl. Math. 188
(2015) 16-24.], and we present some results which correct their
mistakes. We also give the second-minimum and third-minimum degree
resistance distances among graphs in $Cact(n;t)$, and characterize
the corresponding extremal graphs as well.
\end{abstract}



\section{ Introduction}

The graphs considered in this paper are finite, loopless, and
contain no multiple edges. Given a graph $G$, let $V(G)$ and
$E(G)$ be, respectively, its vertex and edge sets. The ordinary
distance $d(u,v)= d_G(u,v)$ between the vertices $u$ and $v$ of
the graph $G$ is the length of the shortest path between $u$ and
$v$~\cite{Bondy1976}. For other undefined notations and
terminology from graph theory, the readers are referred
to~\cite{Bondy1976}.

The Wiener index $W(G)$is the sum of ordinary distances between
all pairs of vertices, that is, $W(G)=\sum_{\{u,v\}\subseteq
V(G)}d(u,v).$ It is the oldest and one of the most thoroughly
studied distance-based graph invariant.
 A modified version of the
Wiener index is the degree distance defined as
$D(G)=\sum_{\{u,v\}\subseteq V(G)}[d(u)+d(v)]d(u,v)$, where $d(u)=
d_G(u)$ is the degree of the vertex $u$ of the graph $G.$

In 1993, Klein and Randi\'c~\cite{Klein1993} introduced a new
distance function named resistance distance, based on the theory
of electrical networks. They viewed $G$ as an electric network $N$
by replacing each edge of $G$ with a unit resistor. The resistance
distance between the vertices $u$ and $v$ of the graph $G$,
denoted by $R(u,v)$, is then defined to be the effective
resistance between the nodes $u$ and $v$ in $N$. If the ordinary
distance is replaced by resistance distance in the expression for
the Wiener index, one arrives at the Kirchhoff
index~\cite{Klein1993,DSTG2015} \[Kf(G)=\sum\limits_{\{ u,v\}
\subseteq V(G)} {R(u,v)}
\]which has been widely studied~\cite{Gao2011,Gao2012,Feng2012,AFeng2012,Feng2014,FengY2014,Bu2014}. In 1996, Gutman and
Mohar~\cite{Gutman1996} obtained the famous result by which a
relationship is established between the Kirchhoff index and the
Laplacian spectrum: $Kf(G)=n\sum^{n-1}_{i=1}\frac{1}{\mu_i},$
where $\mu_1\geq \mu_2\geq \cdots \geq \mu_n=0 $ are the
eigenvalues of the Laplacian matrix of a connected graph $G$ with
$n$ vertices. For more details on the Laplacian matrix, the
readers are referred to~\cite{Liu2015,LiuP2015}. Bapat et al. has
provided a simple method for computing the resistance distance
in~\cite{Bapat2003}.
Palacios~\cite{PalaciosJ2010,Palacios2001,Palacios2004,Palacios2011,Palacios2010,Palacios2015}
studied the resistance distance and the Kirchhoff indices of
connected undirected graphs with probability methods. E. Bendito
et al.~\cite{Bendito2008} formulated the Kirchhoff index based on
discrete potential theory. M. Bianchi et al. obtained the upper
and lower bounds for the Kirchhoff index $Kf(G)$ of an arbitrary
simple connected graph $G$ by using a majorization
technique~\cite{Bianchi2013}. Besides, the Kirchhoff indices of
some lattices are investigated
in~\cite{LY2015,LiuYan,YLZ,ZZ,Liu2014}. Similarly, if the ordinary
distance is replaced by resistance distance in the expression for
the degree distance, then one arrives at the degree resistance
distance\[{D_R}(G) = \sum\limits_{\{ u,v\}  \subseteq V(G)} {[d(u)
+ d(v)]R(u,v)} .\]  Palacios~\cite{Palacios2013} named the same
graph invariant ¡®¡®additive degree-Kirchhoff index¡¯¡¯.

Tomescu~\cite{Tomescu2008} determined the unicyclic and bicyclic
graphs with minimum degree distance. In~\cite{Tomescu2009}, the
author investigated the properties of connected graphs having
minimum degree distance.
 Bianchi et al.~\cite{Bianchi2013} gave some upper and lower bounds for $D_R$ whose expressions do not depend on the resistance distances.
 Yang and Klein gave formulae for the degree resistance distance of the subdivisions and triangulations of graphs~\cite{YangK2015}.
 For more work on $Kf(G)$, the readers are referred
to~\cite{XuG2014,Yang2014,YangJ2008,YangZhang2008,YangK2013,Li2015}.

A graph $G$ is called a cactus if each block of $G$ is either an
edge or a cycle. Denote by $Cact(n;t)$ the set of cacti possessing
$n$ vertices and $t$ cycles~\cite{Hua2010,liu2007}. In this paper,
we determine the minimum degree resistance distance among graphs
in $Cact(n;t)$ and characterize the corresponding extremal graphs.

\section{Preliminaries}

Let ${R_G}(u,v)$ denote the resistance distance between $u$ and$v$ in the graph $G$.
Recall that ${R_G}(u,v) = {R_G}(v,u)$ and ${R_G}(u,v) \ge 0$ with equality if and only if $u = v$.

For a vertex $u$ in $G$,we define
\begin{center}
  $K{f_v}(G) = \sum\limits_{u \in G} {{R_G}(u,v)}$ ~~~and~~~ ${D_v}(G) = \sum\limits_{u \in G} {{d_G}(u){R_G}(u,v)}$.
\end{center}
In what follows, for the sake of conciseness, instead of $u \in
V(G)$ we write $u \in G $. By the definition of ${D_v}(G)$, we
also have\[{D_R}(G) = \sum\limits_{v \in G} {{d_G}(v)}
\sum\limits_{u \in G} {{R_G}(u,v).} \]
\begin{lemma} (\cite{Feng2012}).
  Let $G$ be a connected graph with a cut-vertex $v$ such that $G_1$ and $G_2$ are two connected
  subgraphs of $G$ having $v$ as the only common vertex and $V({G_1}) \cup V({G_2}) = V(G)$.

  Let ${n_1} = \left| {V({G_1})} \right|, {n_2} = \left| {V({G_2})} \right|$,
  ${m_1} = \left| {E({G_1})} \right|$,
  ${m_2} = \left| {E({G_2})} \right|$.
  Then \[{D_R}(G) = {D_R}({G_1}) + {D_R}({G_2}) + 2{m_2}K{f_v}({G_1}) + 2{m_1}K{f_v}({G_2})
   + ({n_2} - 1){D_v}({G_1}) + ({n_1} - 1){D_v}({G_2}).\]
\end{lemma}

\begin{definition}(\cite{DSTG2015}).
  Let $v$ be a vertex of degree $p+1$ in a graph $G$,
  such that $vv_1, vv_2, \ldots, vv_p$ are pendent edges incident with $v$,
   and $u$ is the neighbor of $v$ distinct from $v_1,v_2,\ldots ,v_p$.
   We form a graph ${G^{'}} = \sigma (G,v)$ by deleting the edges $vv_1,vv_2,\ldots,vv_p$
   and adding new edges $ uv_1,uv_2,\ldots,uv_p $.
   We say that ${G^{'}}$ is a $\sigma$-transform of the graph $G$.
\end{definition}

\qquad\quad \includegraphics[width=12cm]{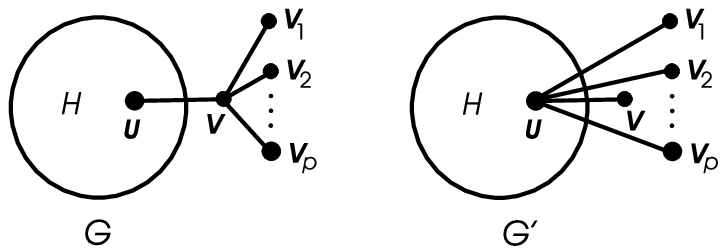}

\qquad\qquad\qquad\qquad\qquad\qquad\qquad The $\sigma$ -transform at $v$

\begin{lemma}(\cite{Feng2012}).
 Let ${G^{'}} = \sigma (G,v)$ be a $\sigma$ -transform of the graph G, ${d_G}(u) \ge 1$.
  Then ${D_R}(G) \ge {D_R}({G^{'}})$. Equality holds if and only if G is a star with v as its center.
\end{lemma}

\begin{lemma} (\cite{Feng2012}).
 Let $C_k$ be the cycle of size k,
 and $v \in C_k$.
 Then,
 $Kf({C_k}) = \frac{{{k^3} - k}}{{12}},
 {D_R}({C_k}) = \frac{{{k^3} - k}}{3},
 K{f_v}({C_k}) = \frac{{{k^2} - 1}}{6}$ and ${D_v}({C_k}) = \frac{{{k^2} - 1}}{3}$.
\end{lemma}


\begin{definition}(\cite{DSTG2015}).
  Let $G \in Cact(n;t), t\geq2$.
  A cycle $C$ of $G$ is said to be an end cycle if there is a unique vertex $v$ in $C$ which is adjacent to a vertex in $V(G)\setminus V(C)$.
  This unique vertex $v$ in $C$ is called the anchor of $C$.
\end{definition}

\begin{lemma} (\cite{DSTG2015}).
Let $G\in Cact(n;t),$ $t\geq2$, be a cactus without cut edges. Let
$C$ be an end cycle of $G$ and $v$ be its anchor. Let $u$ be a
vertex of $C$ different from $v$. The graphs $G_1$ and $G_2$ are
constructed by adding $r$ pendent edges to the vertices $u$ and
$v$, respectively. Then ${G_R}(G_1)>{G_R}(G_2)$.
\end{lemma}

\qquad\qquad \includegraphics[width=12cm]{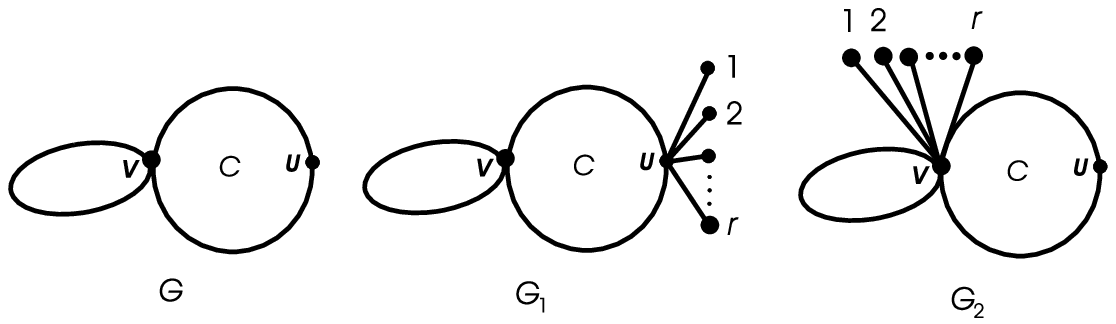}

\section{ Some errors in~\cite{DSTG2015} and corrections}

In~\cite{DSTG2015}, J. Du, G. Su, J. Tu, I. Gutman proved that
$G^0(n;t)$ is the unique element of $Cact(n; t)$, $t\geq1$, having
minimum degree resistance distance. Unfortunately, there are some
computational errors in the process of the proof. We shall list
the errors in~\cite{DSTG2015} as Errors 3.1, 3.2 below.

\qquad\qquad\qquad\qquad\qquad \includegraphics[width=8cm]{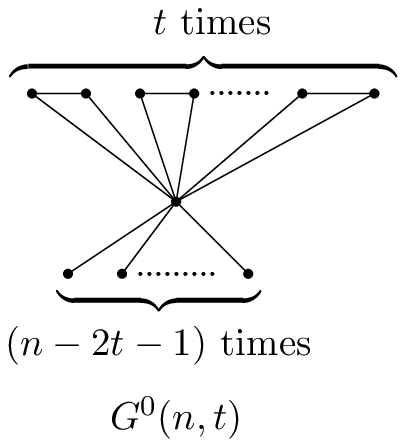}

\vspace{5pt} \noindent {\bf  Error 3.1} (\textbf {Lemma 7
in~\cite{DSTG2015}})

${D_R}({C_h})-{D_R}(S)=\frac{{{h^2} - 8h + 3}}{3}$ and
$|V(H)|-1=n-h-1.$

\vspace{5pt} \noindent {\bf  Counterexample 1}

If $h=4$, according to the Lemma 7 in ~\cite{DSTG2015}, the result
is $-\frac{13}{3}$ and $n-5$. In fact, the correct result is
$-\frac{10}{3}$ and $n-4$, which arrives at a contradiction.

\vspace{5pt} \noindent {\bf  {Correction of Lemma 7 in ~\cite{DSTG2015}}}

 Let $G = (V,E)$ be a graph belonging to
$Cact(n;t)$, $t\geq3$. Let $C_h$ be a cycle with $h\geq4$
vertices, contained in G. Let there be a unique vertex $u$ in
$C_h$ which is adjacent to a vertex in $V(G)\setminus V(C)$.
Assuming that $uv, vw\in E(C)$, construct a new graph
$G^{*}=G-vw+uw$ as shown in the following figure. Then,
${D_R}(G)>{D_R}(G^{*})$.

\qquad\qquad\quad \includegraphics[width=10cm]{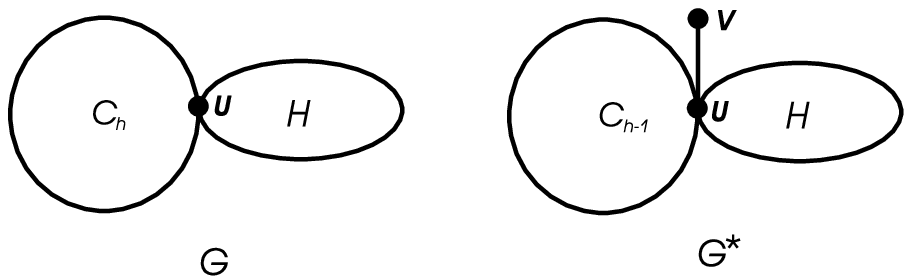}

Let $S$ be the graph obtained by attaching to the vertex $u$ of
$C_{h-1}$ the pendent vertex $v$.
${D_R}({C_h})-{D_R}(S)=\frac{{{h^2} - 8h + 6}}{3}$ and
$|V(H)|-1=n-h$.

Using Lemma 1, we have
\[{D_R}(G) = {D_R}({C_h})+{D_R}(H) + 2|E(H)|K{f_u}({C_h}) + 2hK{f_u}(H) + (|V(H) - 1|){D_u}({C_h}) + (h - 1){D_u}(H),\]
\[{D_R}({G^*}) = {D_R}(S) + {D_R}(H) + 2|E(H)|K{f_u}(S) + 2hK{f_u}(H) + (|V(H) - 1|){D_u}(S) + (h - 1){D_u}(H).\]

Then
\begin{equation}
  \begin{split}
    &{D_R}(G) - {D_R}({G^*})\\
     &= {D_R}({C_h}) - {D_R}(S) + 2(n + t - 1 - h)[K{f_u}({C_h}) - K{f_u}(S)] + (n - h)[{D_u}({C_h}) - {D_u}(S)]\\
    &=\frac{{{h^2} - 8h + 6}}{3} + 2(n + t - 1 - h)\frac{{2h - 7}}{6} + (n - h)\frac{{2h - 4}}{3}\\
    &=\frac{{{h^2} - 8h + 6}}{3}+(n  - 1 - h)\frac{{4h - 11}}{3}+t\frac{{2h - 7}}{3}+\frac{{2h - 4}}{3}\\
    &\geq\frac{{{h^2} - 19}}{3}+(n  - 1 - h)\frac{{4h - 11}}{3} ~~~(by \quad t\geq3).\nonumber
  \end{split}
\end{equation}

If $h=4$, then ${D_R}(G) -
{D_R}({G^*})\geq\frac{5}{3}n-\frac{28}{3}>0$.

If $h\geq 5$, then ${D_R}(G) - {D_R}({G^*})>(n  - 1 - h)\frac{{4h -
11}}{3}>0$.

This completes the proof. \qed

\vspace{5pt} \noindent {\bf  Error 3.2} (\textbf {Theorem 1  in ~\cite{DSTG2015}})

\[{D_R}({G^0}(n,t)) =  - \frac{4}{3}{t^2} - (\frac{8}{3}n - 6)t + 3{n^2} - 7n + 4.\]

\vspace{5pt} \noindent {\bf  Counterexample 2}

If $n=5,t=1$, according to the  Theorem 1 in~\cite{DSTG2015}, the
result is 50. In fact, the correct result is $44\frac{2}{3}$,
which also arrives at a contradiction.

\vspace{5pt} \noindent {\bf  Correction of Error 3.2 }

It is obvious that the ${D^0}(n,t)$ consists of $n$ $C_3$ and an
$S_{n-2t}$, in which $n$ $C_3$ and an $S_{n-2t}$ have a common
vertex $v_1$. Using Lemma 1, we have
\begin{equation}
  \begin{split}
{D_R}({G^0}(n,t)) &= t{D_R}({C_3}) + {D_R}({S_{n - 2t}}) + 2t(n + t - 4)K{f_{{v_1}}}({C_3}) + 6tK{f_{{v_1}}}({S_{n - 2t}})\\
 &~~~+ t(n - 3){D_{{v_1}}}({C_3}) + 2t{D_{{v_1}}}({S_{n - 2t}})\\
 &=8t + (n - 2t)(n - 2t - 1) + 2(n - 2t - 1)(n - 2t - 2) + \frac{8}{3}t(n + t - 4) \\
 &~~~+ 6t(n - 2t - 1) + \frac{8}{3}t(n - 3) + 2t(n - 2t - 1)\\
 &= - \frac{4}{3}{t^2} + (\frac{4}{3}n - \frac{{14}}{3})t + 3{n^2} - 7n + 4.\nonumber
  \end{split}
\end{equation}
\qed

In the following we shall consider the cacti with the second and
the third-minimum degree resistance distances.

\section{The second-minimum degree resistance distance}

By Lemmas 2.2, 2.4 and Theorem 7 in~\cite{DSTG2015}, one can
conclude that $G$ which has the second-minimum degree resistance
distance in $Cact(n;t)$ must be one of the graphs $G_3$, $G_4$,
and $G_5$ as shown in the Figure 1.

\includegraphics[width=14cm]{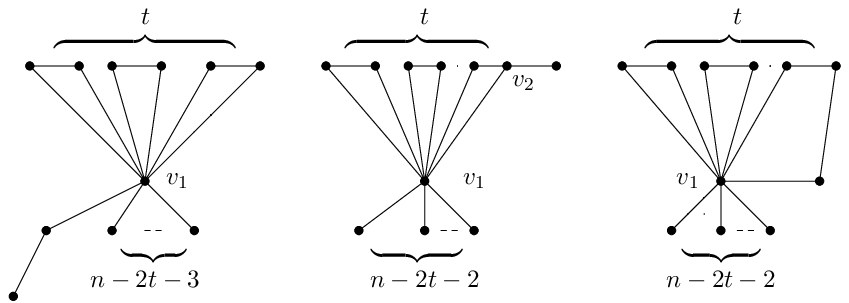}

\qquad\qquad\qquad\quad\quad\qquad$G_3$\qquad\quad\qquad\qquad\quad$G_4$\qquad\qquad\qquad\quad\qquad$G_5$
\begin{center}
  \bf {Figure 1}
\end{center}

\begin{theorem}
  Among all graphs in $Cact(n,t)$ with $n\geq7$ and $t\geq1$,
  the cactus with the second-minimum degree resistance distance is $G_5$.
\end{theorem}
\textbf{Proof}. (i): Let $H_1$ denote the common subgraph of $G_3$
and ${G^0}(n,t)$.
Thus, we can view graphs $G_3$ and ${G^0}(n,t)$
as the graphs depicted in Figure 2.

\qquad\qquad\includegraphics[width=12cm]{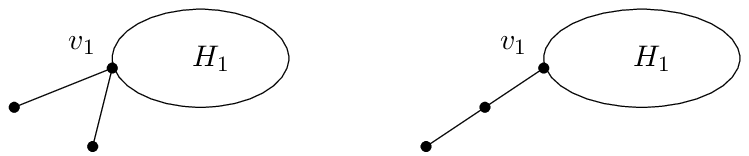}

\qquad\qquad\qquad\qquad\quad\qquad\qquad${G^0}(n,t)$\qquad\qquad\qquad\qquad\qquad\qquad \qquad\quad${G_3}$
\begin{center}
  \bf {Figure 2}
\end{center}

Using Lemma 1, we have
$${D_R}({G^0}(n,t))  = {D_R}({H_1}) + {D_R}({S_3}) + 4K{f_{{v_1}}}({H_1}) + 2(n + t - 3)K{f_{{v_1}}}({S_3}) + 2{D_{{v_1}}}({H_1}) + (n - 3){D_{{v_1}}}({S_3}),$$
$${D_R}({G_3}) = {D_R}({H_1}) + {D_R}({P_3}) + 4K{f_{{v_1}}}({H_1}) + 2(n + t - 3)K{f_{{v_1}}}({P_3}) + 2{D_{{v_1}}}({H_1}) + (n - 3){D_{{v_1}}}({P_3}).$$

Here\quad $K{f_{{v_1}}}({S_3})=2, K{f_{{v_1}}}({P_3})=3,{D_{{v_1}}}({S_3})=2,{D_{{v_1}}}({P_3})=4.$

Therefore,
\begin{equation}
\begin{split}
{D_R}({G_3})-{D_R}({G^0}(n,t)) &=2(n + t - 3)(K{f_{{v_1}}}({P_3})-K{f_{{v_1}}}({S_3}))+(n - 3)({D_{{v_1}}}({P_3})-{D_{{v_1}}}({S_3}))\\
&= 2(n + t - 2) + 2(n - 3)\\
 &= 4n + 2t - 12.\nonumber
\end{split}
\end{equation}

(ii): Let $H_2$ denote the common subgraph of $G_3$ and
${G^0}(n,t)$. Thus, we can view graphs $G_3$ and ${G^0}(n,t)$ as
the graphs depicted in Figure 3.

\qquad\qquad\includegraphics[width=12cm]{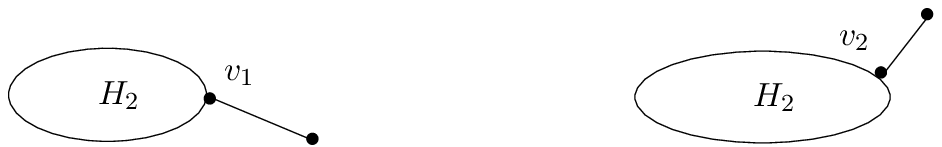}

\qquad\qquad\qquad\qquad\qquad\qquad${G^0}(n,t)$\qquad\qquad\qquad\qquad\qquad\qquad \qquad\quad${G_4}$
\begin{center}
  \bf {Figure 3}
\end{center}

Using Lemma 1, we have
 $${D_R}({G^0}(n,t))  = {D_R}({H_2}) + {D_R}({P_2}) + 2K{f_{{v_1}}}({H_2})
 + 2(n + t - 2)K{f_{{v_1}}}({P_2}) + {D_{{v_1}}}({H_2}) + (n - 2){D_{{v_1}}}({P_2}),$$
 $${D_R}({G_4})= {D_R}({H_2}) + {D_R}({P_2}) + 2K{f_{{v_2}}}({H_2})
  + 2(n + t - 2)K{f_{{v_2}}}({P_2}) + {D_{{v_2}}}({H_2}) + (n - 2){D_{{v_2}}}({P_2}).$$

 Here
 $$K{f_{{v_1}}}({H_2}) = n - \frac{2}{3}t - 2,
 {Kf_{{v_2}}}({H_2}) = \frac{5}{3}n - \frac{2}{3}t - \frac{{14}}{3},$$
 $${D_{{v_1}}}({H_2}) = n + \frac{2}{3}t - 2,
 {D_{{v_2}}}({H_2}) = \frac{7}{3}n + 2t - \frac{{26}}{3}.$$

 Therefore,
 \begin{equation}
  \begin{split}
 {D_R}({G_4})-{D_R}({G^0}(n,t))& =2(K{f_{{v_2}}}({H_2})-K{f_{{v_1}}}({H_2}))+{D_{{v_2}}}({H_2})-{D_{{v_1}}}({H_2})\\
 & =2(\frac{2}{3}n - \frac{8}{3})+(\frac{4}{3}n + \frac{4}{3}t - \frac{{26}}{3})\\
 &=\frac{8}{3}n + \frac{4}{3}t - 12.\nonumber
\end{split}
\end{equation}

(iii): Let $H_2$ denote the common subgraph of $G_5$ and
${G^0}(n,t)$. Thus, we can represent these graphs as follows
 in Figure 4.

\qquad\qquad\includegraphics[width=4cm]{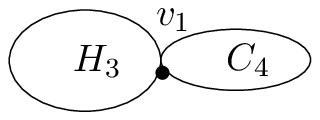}\qquad\qquad\qquad\qquad\includegraphics[width=4cm]{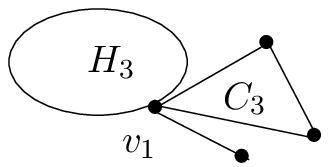}

\qquad\qquad\qquad\qquad\quad${G^0}(n,t)$\qquad\qquad\qquad\qquad\qquad\qquad \qquad\quad${G_5}$
\begin{center}
  \bf{Figure 4}
\end{center}

Using Lemma 1, we have
$${D_R}({G^0}(n,t))  = {D_R}({H_3}) + {D_R}(S_4^3) + 8K{f_{{v_1}}}({H_3}) + 2(n + t - 5)K{f_{{v_1}}}({S_4^3}) + 3{D_{{v_1}}}({H_3}) + (n - 4){D_{{v_1}}}({S_4^3}),$$
$${D_R}(G_5)  = {D_R}({H_3}) + {D_R}(C_4) + 8K{f_{{v_1}}}({H_3}) + 2(n + t - 5)K{f_{{v_1}}}({C_4}) + 3{D_{{v_1}}}({H_3}) + (n - 4){D_{{v_1}}}({C_4}).$$

Here\quad
${D_R}({C_4}) = \frac{{70}}{3},{D_R}(S_4^3) = 20,
K{f_{{v_1}}}({C_4}) = \frac{7}{3},K{f_{{v_1}}}(S_4^3) = \frac{5}{2},
{D_{{v_1}}}({C_4}) = \frac{{11}}{3},{D_{{v_1}}}(S_4^3) = 5$.

Therefore,
 \begin{equation}
  \begin{split}
 {D_R}({G_5})-{D_R}({G^0}(n,t))& ={D_R}(C_4)-{D_R}(S_4^3)+2(n+t-5)(K{f_{{v_1}}}({C_4})-K{f_{{v_1}}}({S_4^3}))\\
 &\quad+(n-4)({D_{{v_1}}}({C_4})-{D_{{v_1}}}(S_4^3))\\
 &=-\frac{10}{3}+\frac{1}{3}(n+t-5)+\frac{4}{3}(n-4)\\
 &=\frac{5}{3}n+\frac{t}{3}-\frac{31}{3}.\nonumber
\end{split}
\end{equation}

By the above expressions for the degree resistance distances of
$G_3$, $G_4$ and $G_5$, we immediately have the desired result.
\qed

From Theorem 4.1 we immediately have the following result.

\begin{corollary}
  For a graph G, not isomorphic to ${G^0}(n,t)$,
  in $Cact(n,t)$ with $n\geq7$ and $t\geq1$,
  it holds that ${D_R}(G) \ge  - \frac{4}{3}{t^2} + (\frac{4}{3}n - \frac{{13}}{3})t + 3{n^2} - \frac{{16}}{3}n - \frac{{19}}{3}$,
  with equality if and only if $G\cong{G_5}$.
\end{corollary}

\section{The third-minimum degree resistance distance}

By the same reasonings as was used in Theorem 4.1, we conclude
that the possible candidates having the third-minimum degree
resistance distance must come from one of $G_4$, $G_6-G_{10}$.

\begin{theorem}
  Among all graphs in $Cact(n,t)$ with $n\geq25$ and $t\geq1$,
  the cactus with the third-minimum degree resistance distance is $G_4$.
\end{theorem}
\textbf{Proof}.\  By above discussions, we need only to determine
the minimum cardinality among $D_R(G_4)$, $D_R(G_6)$, $D_R(G_7)$,
$D_R(G_8)$, $D_R(G_9)$ and $D_R(G_{10})$.

\qquad\quad\qquad\includegraphics[width=4cm]{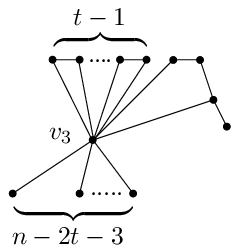}\includegraphics[width=4cm]{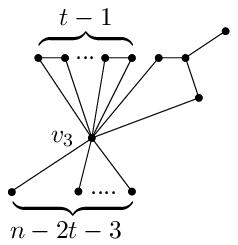}\includegraphics[width=4.3cm]{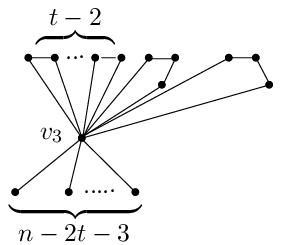}

\qquad\qquad\;\qquad\quad\qquad$G_6$\qquad\quad\qquad\qquad\quad$G_7$\qquad\qquad\qquad\quad\qquad$G_8$

\qquad\quad\qquad\quad\includegraphics[width=4.6cm]{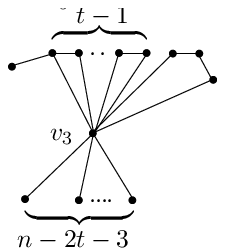}\qquad\qquad\quad\includegraphics[width=3.35cm]{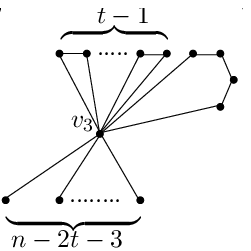}

\qquad\qquad\qquad\;\qquad\quad\qquad$G_9$\qquad\quad\qquad\quad\quad\qquad\qquad\quad$G_{10}$
\begin{center}
  \bf{Figure 5}
\end{center}

Let $H_4$ denote the common subgraph of $G_4$, $G_6$ and $G_7$.
Thus, we can view graphs $G_4,$ $G_6$ and $G_7$ as the graphs
depicted in Figure 6.

\includegraphics[width=4cm]{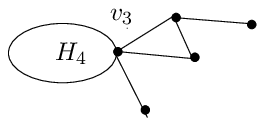}\quad\includegraphics[width=4cm]{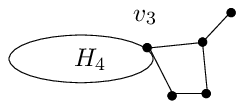}\quad\includegraphics[width=4cm]{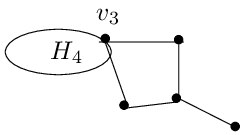}\quad\includegraphics[width=2.2cm]{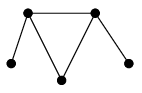}

\qquad\quad$G_4$\qquad\qquad\qquad\quad\qquad\qquad\quad$G_6$\qquad\qquad\qquad\qquad\qquad$G_7$\qquad\qquad\qquad\qquad$G_0$
\begin{center}
  \bf{Figure 6}
\end{center}

Using Lemma 1, we have
$${D_R}({G_4})= {D_R}({H_4}) + {D_R}({G_0}) + 10K{f_{{v_3}}}({H_4}) + 2(n + t - 6)K{f_{{v_3}}}({G_0}) + 4{D_{{v_3}}}({H_4}) + (n - 5){D_{{v_3}}}({G_0}),$$
$${D_R}({G_6})= {D_R}({H_4}) + {D_R}({S_5^4}) + 10K{f_{{v_3}}}({H_4})
 + 2(n + t - 6)K{f_{{v_3}}}({S_5^4}) + 4{D_{{v_3}}}({H_4}) + (n - 5){D_{{v_3}}}({S_5^4}).$$

Here\quad ${D_R}({G_0})=\frac{142}{3},{D_R}({S_5^4})=43,K{f_{{v_3}}}({G_0})=4,
K{f_{{v_3}}}({S_5^4})=\frac{17}{4},
{D_{{v_3}}}({G_0})=6,{D_{{v_3}}}({S_5^4})=\frac{15}{2}.$

Therefore,
\begin{equation}
  \begin{split}
 {D_R}({G_6})-{D_R}({G_4})& ={D_R}({S_5^4})-{D_R}({G_0})+2(n+t-6)(K{f_{{v_3}}}({S_5^4})-K{f_{{v_3}}}({G_0}))\\
 &\quad+(n-5)({D_{{v_3}}}({S_5^4})-{D_{{v_3}}}({G_0}))\\
 &=-\frac{13}{3}+\frac{1}{2}(n+t-6)+\frac{3}{2}(n-5)\\
 &=2n+\frac{t}{2}-\frac{89}{6}>0.\nonumber
\end{split}
\end{equation}

Using Lemma 1, we have
$${D_R}({G_7})= {D_R}({H_4}) + {D_R}({S_5^4}) + 10K{f_{{v_3}}}({H_4}) + 2(n + t - 6)K{f_{{v_3}}}({S_5^4})
+ 4{D_{{v_3}}}({H_4}) + (n - 5){D_{{v_3}}}({S_5^4}).$$

Here\quad $K{f_{{v_3}}}({S_5^4})=\frac{9}{2},{D_{{v_3}}}({S_5^4})=8$.

Therefore,
\begin{equation}
  \begin{split}
  {D_R}({G_7})-{D_R}({G_6})&=\frac{1}{2}(n+t-6)+\frac{1}{2}(n-5)\\\nonumber
  &=n+\frac{1}{2}t-\frac{11}{2}>0.
  \end{split}
\end{equation}

Then\quad ${D_R}({G_7})>{D_R}({G_6})>{D_R}({G_4}).$

Similar to the relationship between ${D_R}({G_5})$ and
${D_R}({G^0}(n,t))$, we have

\begin{equation}
  \begin{split}
  {D_R}({G_8})-{D_R}({G_5})&={D_R}({G_5})-{D_R}({G^0}(n,t))\\\nonumber
  &=\frac{5}{3}n+\frac{t}{3}-\frac{31}{3}.
  \end{split}
\end{equation}
$${D_R}({G_8})-{D_R}({G^0}(n,t))=\frac{10}{3}n+\frac{2}{3}t-\frac{62}{3}.$$

Therefore,
$${D_R}({G_8})-{D_R}({G_4})=\frac{2}{3}n-\frac{2}{3}t-\frac{26}{3}.$$

Because of $t\leq\frac{n-1}{2}$, when $n\geq25$, ${D_R}({G_8})-{D_R}({G_4})>0$.

Similar to the relationship between ${D_R}({G_5})$ and
${D_R}({G^0}(n,t))$, we have
$${D_R}({G_9})-{D_R}({G_4})={D_R}({G_5})-{D_R}({G^0}(n,t))=\frac{5}{3}n+\frac{t}{3}-\frac{31}{3}.$$
$${D_R}({G_9})-{D_R}({G^0}(n,t))=\frac{13}{3}n+\frac{5}{3}t-\frac{67}{3}.$$

Therefore,
$${D_R}({G_9})-{D_R}({G_8})=n+t-\frac{5}{3}>0.$$

Then\quad${D_R}({G_9})>{D_R}({G_8})>{D_R}({G_4}).$

Similar to the method of ${D_R}({G_5})-{D_R}({G^0}(n,t)),$ we have
\begin{equation}
  \begin{split}
{D_R}({G_{10}})-{D_R}({G_5})& ={D_R}(C_5)-{D_R}(S_5^4)+2(n+t-6)(K{f_{{v_1}}}({C_5})-K{f_{{v_1}}}({S_5^4}))\\
&\quad+(n-5)({D_{{v_1}}}({C_5})-{D_{{v_1}}}(S_5^4))\\
&=-3+(n+t-6)+2(n-5)\\
&=3n+t-19>0.\nonumber
\end{split}
\end{equation}

Then\quad ${D_R}({G_{10}})>{D_R}({G_5})>{D_R}({G_4}).$

By the
above several inequalities, we immediately have the desired
result.\qed

From Theorem 5.1 we immediately have the following result.

\begin{corollary}
    For a graph G, not isomorphic to ${{G^0}(n,t),G_5}$, in $Cact(n,t)$ with $n\geq25$ and $t\geq1$,
    it holds that ${D_R}(G) \ge  - \frac{4}{3}{t^2} + (\frac{4}{3}n - \frac{{10}}{3})t + 3{n^2} - \frac{{13}}{3}n - 8$,
    with equality if and only if $G\cong{G_4}$.
\end{corollary}

 \vspace{5pt} \noindent
{\bf Acknowledgments}\

The work of J. B. Liu was partly supported by the Natural Science
Foundation of Anhui Province of China under Grant No. KJ2013B105
and the National Science Foundation of China under Grant Nos.
11471016 and 11401004; The work of X. F. Pan was partly supported
by the National Science Foundation of China under Grant Nos.
10901001, 11171097, and 11371028.

\end{document}